\documentclass{article}
\usepackage{amssymb,amsfonts,amsmath}

\title{The Maximum or Minimum Number \\ of Rational Points on 
 Curves \\ of Genus Three over Finite Fields}
\author{Kristin Lauter\\
with Appendix by Jean-Pierre Serre}
\date{}

\def\C{\mathbb{C}}
\def\F{\mathbb{F}}

\def\Z{\mathbb{Z}}
\def\Q{\mathbb{Q}}

\newcommand{\End}{\mbox{End}}
\newcommand{\Ab}{\mbox{Ab}}
\newcommand{\Mod}{\mbox{Mod}}
\newcommand{\Hom}{\mbox{Hom}}
\newcommand{\rang}{\mbox{rang}}
\newcommand{\Cl}{\mbox{Cl}}
\newcommand{\Pic}{\mbox{Pic}}
\newcommand{\Ker}{\mbox{Ker}}
\newcommand{\Aut}{\mbox{Aut}}

\newtheorem{thm}{Theorem}
\newtheorem{lem}{Lemma}
\newtheorem{Def}{Definition}
\newtheorem{cor}{Corollary}

\newtheorem{prop}{Proposition}

\begin{document}
\maketitle


\noindent
\begin{abstract}                
We show that for all finite fields $\F_q$, there exists a
curve $C$ over $\F_q$ of genus $3$ such that the number of
rational points on $C$ is within $3$ of the Serre-Weil
upper or lower bound.  For some $q$, we also obtain
improvements on the upper bound for the number of rational
points on a genus $3$ curve over $\F_q$.
\end{abstract}

\section{Introduction}

More than half a century ago, Andr\'e Weil proved a formula for
the number of rational points, $N(C)$, on a smooth projective algebraic
curve $C$ of genus $g$ over a finite field $\F_q$.  This formula,
along with his proof of what is referred to as the Riemann hypothesis
for curves, provides upper (resp. lower) bounds on the maximum
(resp. minimum) number of rational points possible  
$$ q+1-2g\sqrt{q} \le N \le q+1+2g\sqrt{q}.$$
There are many cases in which the Weil upper and lower bounds cannot
be attained.  Some are trivial: for example when the bound is 
not an integer.  Also, when the field size, $q$, is small
with respect to the genus, $g$, the lower bound will be negative
and thus cannot be attained.
In \cite{Se}, Serre made a non-trivial improvement to the Weil bound 
(which we will refer to hereafter as the Serre-Weil bound):
$$ q+1-gm \le  N \le q+1+gm, \quad m = [2\sqrt{q}],$$
and introduced the explicit formulae method to provide 
better bounds for large genus.  
Since then there has been considerable interest
in determining the actual maximum and minimum.  
(cf. \cite{Se}, \cite{Serre82}, \cite{Serre83}, \cite{Serre-Harvard}, 
\cite{Stark}, \cite{Sc} \cite{Stohr}, \cite{Geer}, 
\cite{La1}, \cite{La2}, \cite{La3}, \cite{Lacr}, 
\cite{LauterJAG}, \cite{NX1}, \cite{NX2}, \cite{Auer}) 

In the present paper we are concerned with the following question
which was posed in \cite{Serre82}: 
for which genus, $g$, is the difference
between the upper bound and the actual maximum $N_q(g)$
bounded as $q$ varies?  For genus $1$ and any $q$, 
the difference is either or $0$ or $1$ (\cite{Waterhouse}). 
For genus $2$, Serre determined $N_q(2)$ for all $q$, and showed that
the difference from the Serre-Weil bound is always less than or equal to $3$
(\cite{Se});  
for genus $3$, he determined the maximum for $q \le 25$ (\cite{Serre83}).
The present paper is devoted to showing that for genus $3$ and all $q$, 
either the maximum or the minimum is within $3$ of the Serre-Weil 
upper or lower bound.

The techniques involved in the proof of the main theorem include
Serre's theory of hermitian modules as well as 
``glueing'' of polarizations on abelian varieties.  The
theory of hermitian modules is detailed in the Appendix.
This theory provides an equivalence of categories between abelian
varieties over $\F_q$ which are isogenous to a product of copies
of an ordinary elliptic curve, and torsion-free modules
of finite type over a ring which is defined in terms of the Frobenius
of the elliptic curve.  A polarization on the abelian
variety then translates to a hermitian form on the module.
Thus the classification of hermitian forms over rings can be used
to determine the existence or non-existence of the corresponding
polarized abelian variety.
  
To determine whether there is a curve whose number of points is 
close to the Weil bound, we first make a list of the possible zeta
functions for such a curve as in \cite{LauterJAG}.
The numerator of the zeta function of a
curve is given by the characteristic polynomial of Frobenius
acting on the Jacobian of the curve.  From Tate's theorem we
know that the isogeny type of an abelian variety
over a finite field is determined by the characteristic polynomial
of Frobenius.  For each zeta function, we investigate the 
corresponding isogeny type using the equivalence of categories
provided by Serre's theory.  Since not all Jacobians of curves we consider
are isogenous to the product of one elliptic curve with itself,
we are naturally led to consider the glueing of polarizations
on abelian varieties of different isogeny types.  In some cases,
the glueing is possible; in others, it
is not.  When it is not, we obtain improvements
on the upper bounds.  In all cases we are able to conclude  
by making use of the Torelli theorem in dimension $3$ that there 
exists a curve whose number of points is within $3$ of the 
Serre-Weil upper or lower bound.

Serre invoked the theory of hermitian modules to treat the genus $2$ case 
in \cite{Serre82} and \cite{Serre83}, but did not give details.
In \cite{Serre-Harvard}, the full proof of the genus $2$ case is
given, including more explanation of the equivalence of categories. 
More recent work on this subject was done by Everett Howe (\cite{Howe}), 
who translated the notion of principal polarization working with 
Deligne's more general equivalence of categories. The notion
of the glueing of polarizations on abelian varieties has also
been employed by several authors recently for various purposes
(\cite{Kani}, \cite{Poonen}).

The proof of the main theorem of this paper is divided into cases which correspond
naturally to the existence or non-existence of indecomposable hermitian forms
over certain rings of a given discriminant. 
In each case we treat finite fields $\F_q$ with $q$ of a special
form, which leads to a number of interesting diophantine problems. 

The paper is organized as follows:
Section $2$ contains the statement of the main theorems
of the paper.  Section $3$ contains a description
of how the theory of hermitian modules will be used in the
proofs of the main theorems and also recalls the short list
of possible zeta function for curves with the number of points
under consideration.  Section $4$ contains the proofs of the
main theorems.

{\bf Acknowledgements} I would like to thank J-P. Serre for
his letters on hermitian modules which inspired this project,
and for his generous help and time spent working on it.
This paper is based on work supported by a Mentoring Grant
from the Association for Women in Mathematics. 
I would also like to thank Jim Milne and Ren\'e Schoof 
for helpful discussions and Fritz Beukers for the reference
to Ap\'ery's paper. 

\section{Statement of Results}

\noindent
Let $q=p^e$, with $p$ prime, $e \ge 1$.
By a {\it curve} over the finite field $\F_q$, we mean a smooth, projective,
absolutely irreducible curve.  For such a curve, $C$,
let $g=g(C)$ denote the genus, and $N=N(C)$ denote the
number of rational points over $\F_q$.

\vspace{10pt} 
\noindent
{\bf Definition } For fixed $g$ and $q$, let $N_q(g)$ (resp. $M_q(g)$) denote
the maximum (resp. minimum) of $N(C)$ as $C$ runs through all curves of genus
$g$ over $\F_q$.

\vspace{10pt} 
\noindent
Throughout we write $q=p^e$ uniquely in the form $$q=x^2+x+a,$$
where $x$ is the largest integer whose square is less
than or equal to $q$, and $a$ is an integer such that
$$-x \le a \le x.$$  Let $m=[2\sqrt{q}]$, and $d=m^2-4q$.  
Note that
$$ m
=
\left\{
\begin{array}{cl}
 2x & \mbox{if \quad $-x \le a \le 0$,} \\
 2x+1 & \mbox{if \quad $0 < a \le x$,} \\
\end{array}
\right.
$$
and
$$ d
=
\left\{
\begin{array}{cl}
 -4(x+a) & \mbox{if \quad $-x \le a \le 0$,} \\
  1-4a & \mbox{if \quad $0 < a \le x$.} \\
\end{array}
\right.
$$

\begin{thm}  \label{d-3-11} 
Suppose $q=x^2+x+a$, $a=1$ or $a=3$, with 
$a \le x$. Then $$N_q(3) \le q+1+3m-3$$ and
$$M_q(3) \ge q+1-3m+3.$$
Furthermore, there exists a curve $C$ of genus $g(C)=3$ 
over $\F_q$ such that $$|N(C)-(q+1)|=3m-3.$$
\end{thm}

\noindent {\bf Note 2.1}
For example if we write $q=5$ in the form $q=x^2+x+3$, with $x=1$ and $a=3$,
then $x$ is not the greatest integer part of the square root of $q$ and
$a$ does not satisfy the inequality $a \le x$;
so instead we write $q$ in the form $q=x^2+1$, $a=-1$.  
Nor do we write $q=9$ in the form $q=x^2+x+3$, 
since it would not satisfy $a \le x$. 

\begin{thm}   \label{d-4-8}
Suppose $q=x^2+b$, $b=1$ or $b=2$, with $a=b-x$ 
satisfying $-x \le a \le 0$.
Then $N_q(3) \le q+1+3m-2.$
Furthermore, there exists a genus $3$ curve $C$ over $\F_q$
such that $$|N(C)-(q+1)|=3m-2.$$
\end{thm}

\noindent {\bf Note 2.2}
For example $q=3$ is {\it not} of the form $q=x^2+2$ with $2 \le x$;
instead it is of the form $q=x^2+x+1$.

Together, Theorems \ref{d-3-11} and \ref{d-4-8}
can be used to establish the following result.

\begin{thm} \label{all} For any prime power $q=p^e$,
there exists a curve $C$ of genus $g(C)=3$ over $\F_q$
such that $$|N(C)-(q+1)| \ge 3m-3.$$
\end{thm}
In other words, Theorem \ref{all} says that for all $q$,
at least one of the following holds:
$$ i. \quad |q+1+3m-N_q(3)| \le 3,$$
$$ii. \quad |q+1-3m-M_q(3)| \le 3.$$

\section{Background}
\subsection{Hermitian modules}

In the Appendix, Serre gives an equivalence between the
following two categories:
the category of abelian varieties over $\F_q$ which
are isogenous over $\F_q$ to a product of copies of $E$,
where $E$ is an ordinary elliptic curve over $\F_q$;
and the category of torsion-free $R_d$-modules of finite
type, where $R_d=\Z[\pi]$, $\pi$ is the Frobenius of $E$,
and $d=a^2-4q$, where $\#E(\F_q) = q + 1 - a$.
The equivalence holds under the assumption that $d$ is
the discriminant of an imaginary quadratic field. In that
case, $R_d$ is equal to the ring of integers in the
field; thus it is a Dedekind domain, and the modules under
consideration are projective.  If $d$ is not the
discriminant of an imaginary quadratic field,
then $R_d$ is an order in the ring of integers.
Although the equivalence of categories does not
necessarily hold under those conditions, we will still
use the same notation and make use of the functor $S$ 
given in the Appendix.
In most cases there will be no conflict with the notation
from Section $2$ since we will let $a=m$; in all other cases,
we will use the notation $d'$.  The condition that $E$ be
ordinary is equivalent to requiring that $a$ be prime to $p$. 
Throughout the paper we will use the notation $E_m$ to 
denote an elliptic curve with $q+1+m$ points over the
field $\F_q$. 

In Section $5$ of the Appendix, polarizations on abelian
varieties are translated into positive definite hermitian
forms on $R$-modules, and the polarization is principal
if and only if the hermitian form has discriminant $1$.
The Jacobian of a curve has a canonical principal polarization which 
corresponds to the theta divisor.  For an absolutely 
irreducible curve, the theta divisor is irreducible, so
the canonical polarization corresponds to an indecomposable
hermitian module with discriminant $1$.

We will use this correspondence in two directions.
In cases where we can show that there is no
indecomposable hermitian module of discriminant $1$,
we can conclude that no curve of that type exists.
In cases where we find an indecomposable hermitian 
module of discriminant $1$, we can use the theorem
of Torelli in its precise form (see the Appendix to \cite{LauterJAG})
to conclude that a curve of that type (or of the opposite type) exists.

\subsection{Zeta functions}
In each case, we identify the zeta function of the curve
we are searching for.  The zeta function determines the
isogeny type of the Jacobian of the curve.

\begin{Def} A curve has zeta function of type $[x_1,...x_g]$
if  $$x_i = -(\alpha_i + \bar{\alpha_i}), \quad i=1,...,g,$$
where $\{\alpha_i, \bar{\alpha_i} \}$ is the family of $g$ conjugate
pairs of eigenvalues of Frobenius acting on the Jacobian of the curve over
$\F_q$.
\end{Def}

\begin{Def}
A curve $C$ has {\it defect $k$} if $N(C)=q+1+gm-k$, $m=[2\sqrt{q}]$. 
\end{Def}

\noindent {\bf Fact 3.1} (Defect $0$) We recall from \cite{Se} 
that a curve which meets the Serre-Weil bound has zeta function 
of type $[m,m,...,m]$.  When $g=3$, by statement $7.1$ in the Appendix, 
there does not exist a curve of type $[m,m,m]$ if 
$$d=m^2-4q=-3,-4,-8, \text{ or }-11.$$

\noindent {\bf Fact 3.2}  (Defect $1$) We recall from \cite{Serre83} 
that a curve whose number of points is equal to $q+gm$ must have $g \le 2$.

\vspace{10pt}
\noindent {\bf Fact 3.3}  (Defect $2$) We recall from 
\cite{Serre-Harvard} or \cite{LauterJAG} 
that if the fractional part of $2\sqrt{q}$, (which we denote by $\{2\sqrt{q}\}$),
satisfies $\{2\sqrt{q}\} < \sqrt{3} -1$, ($g \ne 4$), then a curve whose number
of points is equal to $q+gm-1$ is of type $[m,m,...,m-2]$.

\section{Proofs of Theorems}

This section is organized as follows.  We will prove the
theorems from Section 2 in the order they were 
stated.  For each possible zeta function to be treated via
the equivalence of categories, we must 
check the condition that the trace be relatively prime
to the characteristic, and deal with special cases where
this fails.

To begin, we notice that except for one special case ($q=2$),
none of the fields in Theorems \ref{d-3-11} and  \ref{d-4-8}
have characteristic $2$.

\noindent {\bf Fact 4.1} If $q=p^e$, $e>1$, then $p \ne 2$ if $q$ is of 
any of the following forms:  
\begin{center}
i. \quad $q=x^2+1$, \\
ii. \quad $q=x^2+2$, \\
iii. \quad $q=x^2+x+1$,\\ 
iv. \quad $q=x^2+x+3$.
\end{center}
{\bf Proof:} In case i, if $q$ were
a power of $2$ then $x$ would be odd, and so
$$x^2 \equiv 1 \pmod{8},$$ which implies  $ q \equiv 2 \pmod{8}$.
Thus $q=2$, $e=1$ is the only possibility. 

In case ii, $x$ would be even, and we would
have $$x^2 \equiv 0 \pmod{4},$$  which implies $ q \equiv 2 \pmod{4}.$

In cases iii. and iv, $x^2+x$ is always even, so $q$ is odd. $\square$

\subsection{Proof of Theorem \ref{d-3-11}}

Write $q=p^e=x^2+x+a$, $a=1$ or $a=3$, with $a \le x$.  Then $m=2x+1$,
and $d=-3$ for $a=1$, and $d=-11$ for $a=3$.

\begin{prop}\label{primeto}
If $q$ is of the form $q=p^e=x^2+x+a$, $a=1$ or $a=3$, $a \le x$, then
$m$ is prime to $p$ unless $q=3$.
\end{prop}
{\bf Proof: }
{\bf a=1.}  Write $q=p^e=x^2+x+1$ with  $m=2x+1$.
Suppose $p$ divides $m$.  Then $$p^e=mx-(x^2-1)$$ 
implies that $p$ divides $(x^2-1)$.  If $p$ divides $(x+1)$, then  $p$ divides
$x$, which is impossible.  So $p$ divides $(x-1)$.  
Thus $p$ also divides $(x+2)$ and  $3$. But
Nagell and Ljunggren (see \cite{Riben}, for example) 
have shown that the only solution to
$$p^e = x^2+x+1, \quad e \ge 3, \quad e \text{  odd,  }$$
is $p=7$, $e=3$, $x=18$. So  no odd power of $3$, $e \ge 3$, is of the
form $x^2+x+1$. So $q=3$ is the only exception. Indeed in that case $m=3$.


\vspace{10pt} 
\noindent
{\bf a=3.} Write $q=p^e=x^2+x+3$ and $m=2x+1$.
Again  suppose $p$ divides $m$.
Then $p^e-3m =x(x-5)$ is divisible by $p$.  We cannot have $p$ divides 
$x$ since then we would also have $p$ divides $(x+1)$.  So in fact, $p$
divides $(x-5)$.  In addition, $2p^e-mx =x+6$ is divisible by $p$.
Thus we must have $p=11$.  However, there are no solutions
to the equation $11^e=x^2+x+3$ since there are no solutions modulo $5$. 
$\square$

\begin{prop}\label{primetom-2}
If $q$ is of the form $q=p^e=x^2+x+a$, $a=1$ or $a=3$, $a \le x$, then
$m-2$ is prime to $p$ unless $q=343$ or possibly $p=5$.
\end{prop}
{\bf Proof: }
{\bf a=1.}  Write $q=p^e=x^2+x+1$ with  $m=2x+1$.
Suppose $p$ divides $m-2$.  Then $$2p^e-(m-2)x=3x+2$$ 
implies that $p$ divides $3x+2$ and $x+3$ and $x-4$.  
So $p=7$. If $e=1$, then $p=7$ does not divide $2x-1$.
If $e$ is odd, $e \ge 3$, then due to the result of
Nagell and Ljunggren cited above, the only solution 
is $p=7$, $e=3$, $x=18$. In that case, $p=7$ divides
$35$, so $q=343$ is an exception.
\vspace{10pt} 
\noindent
{\bf a=3.} Write $q=p^e=x^2+x+3$ and $m=2x+1$.
Again  suppose $p$ divides $m-2$.
Then $$(m-2)^2-4p^e = -8x-11$$ is divisible by $p$,
and so is $8x-4$, which implies that $p$ divides $15$.
If $p=3$, then the only powers of  $p$, $e$ odd, which are of the
form  $x^2 + x + 3$  are  $3$  and  $3^5 = 243$. This can be proved (as Serre
pointed out to me) by Skolem's $\ell$-adic method with $\ell = 5$  
(in \cite{Skinner}, Skinner attempted to give an $\ell$-adic proof with 
$\ell = 11$, but his argument is incomplete).
Neither $q=3$ nor $q=243$ satisfy the divisibility condition ``$3$
divides $2x-1$''.
It would not be necessary to exclude $p=5$ in the hypotheses if
there are no solutions to the equation $5^e=x^2+x+3$ with $e \ge 3$, 
$e$ odd and $5$ dividing $2x-1$.  I have checked with the help of pari
that there are no solutions for $e<1600$.
$\square$

\vspace{10pt}
\subsubsection{q=3}
To prove Theorem \ref{d-3-11}, first consider the case $q=3$.
If $q=3$, then $m=3$, and the explicit formula bound is $N(C) \le 10$.
In fact, a curve $C$ of genus $3$ and defect $3$ over $\F_3$ 
with $N(C)=10=q+1+3m-3$
was given in \cite{Serre-Harvard} with the equation $$y^3-y=x^4-x^2.$$
Its zeta funtion is of type $[m,m,m-3]$.

\subsubsection{\bf $q \ne 3$}  

By Fact 3.1, defect $0$ is not possible because there is
no indecomposable rank $3$ Hermitian module of discriminant $1$ when $d=-3$ or
$d=-11$ (cf. \cite{Hoffmann}). 
By Fact 3.2, defect $1$ is never possible for $g>2$.

\begin{prop}  \label{defect2} If $q$ is of the form $q=p^e=x^2+x+a$,  $a=1$ or
$a=3$ with $a \le x$, then a defect $2$ curve of genus $3$ is of
type $[m,m,m-2]$. 
\end{prop}
{\bf Proof:} By Fact 3.3, it suffices to check that the fractional
part of $2\sqrt{q}$ satisfies $\{2\sqrt{q}\} < \sqrt{3} -1$.
Suppose that $\{2\sqrt{q}\} \ge \sqrt{3} -1.$  We can write
$$\{2\sqrt{q}\} = \{\sqrt{4x^2+4x+4a}\} = \sqrt{4x^2+4x+4a} - (2x+1),$$
so $$\sqrt{4x^2+4x+4a} \ge 2x+\sqrt{3}$$ and then $$4a-3 \ge 4x(\sqrt{3}-1).$$
For $a=1$, this implies $x=0$, which is not possible.  
For $a=3$, this implies $x \le 3$, which does not occur for any $q$. 
(See Note 2.1).
$\square$

\begin{prop} \label{HM2} There is no indecomposable rank $2$, Hermitian module
of discriminant $2$ over $R_d$ when $d=-3$ or $d=-11$.
\end{prop}
{\bf Proof:} 
Since the class number of $R_d$ is $1$, projective $R_d$ modules
are free, so we can express the module $P$ and the Hermitian form $H$
as a matrix whose $ij^{th}$ entry is  $H(e_i,e_j)$ where $\{e_i\}$
is a basis for $P$ over $R_d$.  A change of basis will give
an equivalent matrix. We will work with the matrix in a reduced form.
A rank $2$ Hermitian form in ``reduced'' form can be written as
 \begin{displaymath}
       \left[ \begin{array}{cc}
           \lambda & \bar{\alpha} \\
           \alpha & \mu \\
      \end{array} \right]
    \end{displaymath}   
with $0 < \lambda \le \mu \in \Z$,  $\alpha \in R_d$.  
The condition that $H(x,y)=\overline{H(y,x)}$ implies that $\lambda$, 
$\mu \in \Z$, and we assume that $e_1$ is chosen with $\lambda$
minimal so that $\lambda \le \mu$.  We are interested in Hermitian forms 
which are positive definite, so $\lambda > 0$.  

The proof of the proposition relies on the following lemma:

\begin{lem} If $d$ is square-free, satisfies $d \equiv 1 \pmod{4}$
and $R_d$ has class number one, then we can find a basis for $P$ over
$R_d$ such that the matrix for a Hermitian form in reduced form satisfies
$$\frac{\alpha\bar{\alpha}}{\lambda^2}\le (\frac{|d|+1}{4})^2\frac{1}{|d|}.$$
\end{lem}
Proof: If $\frac{\alpha\bar{\alpha}}{\lambda^2}$ is too large, we can replace
$\alpha$ by $\alpha + \lambda r$ for any $r \in R_d$ by replacing the
basis element $e_2$ by $e_2' = e_2+\bar{r}e_1$.
As a complex number,
$$(\alpha +\lambda r)\overline{(\alpha+\lambda r)}=|(\alpha+\lambda r)|^2,$$
and $|\lambda|=\lambda$, so it suffices to show that $r$ can be chosen so that 
$$|\frac{\alpha}{\lambda}+r|^2 \le (\frac{|d|+1}{4})^2\frac{1}{|d|}.$$
So it is enough to show that for every $z \in \C$, there exists
an element  $r \in R_d$ such that the distance squared from $z$ to $r$
is less than the bound.  A $\Z$- basis for $R_d$ is
$\{1,\frac{1+\sqrt{d}}{2}\}$.  We look for the point in the complex plane
which is furthest from a lattice point, or in other words, the smallest
radius so that circles centered at the lattice points will cover the plane.
It suffices to consider points in the right triangle with vertices
$(0,0),(\frac{1}{2},0),(\frac{1}{2},\frac{\sqrt{|d|}}{2})$, since 
we can then extend the argument by symmetry to the rest of the fundamental
domain. We look for the point $(\frac{1}{2},a)$ in the triangle which is
equidistant from the two lattice points, 
$(0,0)\text{ and }(\frac{1}{2},\frac{\sqrt{|d|}}{2})$.
Setting $$\frac{\sqrt{|d|}}{2} -a = \sqrt{\frac{1}{4}+a^2},$$
we find $$a = \frac{(|d|-1)}{4}\frac{1}{\sqrt{|d|}}.$$  Calculating
the distance squared from this point to the origin we find exactly
the bound stated in the lemma, and this point is the furthest
possible distance away from the closest lattice point. c.q.f.d.

Suppose that $\lambda\mu-\alpha\bar{\alpha}=2$.
The form 
\begin{displaymath}
       \left[ \begin{array}{cc}
           1 & 0 \\
           0 & 2 \\
      \end{array} \right]
    \end{displaymath}  
is decomposable.  If $\lambda=1$, then 
$$\alpha\bar{\alpha} \le (\frac{|d|+1}{4})^2\frac{1}{|d|}
=
\left\{
\begin{array}{cl}
\frac{1}{3}  & \mbox{if $d=-3$} \\
\\
\frac{9}{11} & \mbox{if $d=-11$,} \\
\end{array}
\right.
$$
which implies that $\alpha = \bar{\alpha} = 0$ and thus $\mu=2$.
So it suffices to show $\lambda=1$.  Suppose $\lambda \ge 2$.
For $d=-3$,  $$\mu \le \frac{2+\alpha\bar{\alpha}}{\lambda}
\le \frac{2}{\lambda} + \frac{\lambda}{3},$$ which is less than
$\lambda$ if $\lambda \ge 2$.  This contradicts the fact that
$\lambda \le \mu $.  For $d=-11$,  
$$\lambda^2 \le \lambda\mu = 2+\alpha\bar{\alpha} \le
2+\frac{9}{11}\lambda^2.$$
Thus $\lambda^2 \le 11$, so $\lambda=2 \text{ or }3$.

If $\lambda=2$, then $$\alpha\bar{\alpha} \equiv 0 \pmod{2},$$
and since $$\frac{\alpha\bar{\alpha}}{4} < 1,$$ we have
$\alpha\bar{\alpha} = 2$ or $\alpha\bar{\alpha} = 0$.
But $2$ is not the norm of an element of $R_{-11}$, so
$\alpha\bar{\alpha} = \alpha = 0$. Then we must have
$\mu = 1$, and this contradicts $\lambda \le \mu $.
 
If $\lambda=3$, then $$\alpha\bar{\alpha} \le \frac{9}{11}9 < 8.$$
But then $3\mu = 2 + \alpha\bar{\alpha} \le 9$ implies that
$\mu = 3$ and $\alpha\bar{\alpha} = 7$.  Again this is impossible
since $7$ is not the norm of an element of $R_{-11}$.
$\square$

\vspace{10pt}
\noindent
{\bf Note 4.1.2} Proposition \ref{HM2} can also be deduced from the results
of Otremba (see \cite{Otremba}). In that paper, she proves (via the mass 
formula) that some hermitian forms are ``alone in their genus''; i.e. any 
lattice which is locally isomorphic to the given one is globally isomorphic. 
See especially (for rank 2) p.9 and p.13. These computations imply 
Proposition \ref{HM2}, provided one checks
that every lattice with discriminant $2$ is in the same genus as 
the decomposable one (the {1,2} lattice, in Otremba's notation).

\subsubsection{Glueing criteria}

The following section is translated from \cite{Serreletterglue}.
Suppose that $B$ and $C$ are two polarized abelian varieties
over a perfect field $k$
$$b: B \rightarrow B^* \text{  and  } c: C \rightarrow C^*.$$
Suppose that the polarizations have the same degree, $n$, and that
$n$ is prime to the characteristic, $p$, of $k$ (if the characteristic
is $\ne 0$).  Denote by $N_b$ and $N_c$ the kernels of $b$ and $c$;
these are finite \'etale group schemes of order $n^2$.  We will identify
them with their $\bar{k}$ points. The Galois group $G=\text{Gal}(\bar{k}/k)$
operates on these groups.

Let $\mu$ be the group of roots of unity, written additively.
According to Mumford, [\cite{Mumford}, p.227]
the polarizations $b$ and $c$ define
non-degenerate, alternating bilinear forms on $N_b$ and $N_c$,
with values in $\mu$, and compatible with the action of $G$.
We will denote them by $(x,y) \mapsto \langle x, y \rangle.$
Let $f$ be a map $$f: N_b \rightarrow N_c$$
verifying the following conditions:
\begin{enumerate}
\item  $f$ {\it is an isomorphism of $G$-modules};

\item  $\langle f(x), f(y) \rangle = - \langle x, y \rangle$ \it{for all} 
$(x,y) \in N_b \times N_b.$
\end{enumerate}
To the data $(B,b,C,c,f)$ given above, we can now associate
{\it an abelian variety, $A$, isogenous to $B \times C$,
equipped with a polarization $a$ of degree $1$} as follows:

Let $B \times C$ have the polarization which is the product 
of the polarizations $b$ and $c$.  The kernel of 
$B \times C \rightarrow B^* \times C^*$ is $N_b \times N_c$.
Let $F$ be the subgroup of this kernel which is the graph of the
isomorphism $f$.  Property (2) above shows that $F$ is totally
isotropic; property (1) shows that $F$ is stable by the action of
$G$, i.e. that it is a finite \'etale $k$-subgroup of $B \times C$.
According to Mumford \cite{Mumford}, the polarization on $B \times C$ passes to
the quotient by $F$.  In other words, if we set $A = (B \times C)/F$,
there exists a polarization $a$ on $A$ such that $(b,c)$ factors through $a$
$$B \times C \rightarrow A \rightarrow^a A^* \rightarrow B^* \times C^*.$$
Comparing the degrees, we see that $\deg(a)=1$, so $a$ is principal.

This is how we obtain a polarized abelian variety $A$ by ``glueing''
$(B,b)$ to $(C,c)$ via $f$.

\subsubsection{$q \ne 3$ continued}
The general framework for the glueing of polarized abelian varieties
can be applied in both backwards and forward directions.  The following
theorem works in the backwards direction (``unglueing'').

\begin{thm} \label{notglue2}
Suppose $q=p^e$, $q \ne 3$, $q \ne 343$, $p \ne 5$ and that $q$ is of the form 
$q=x^2+x+a$, $a=1$ or $a=3$ with $a \le x$. Let
$m=[2\sqrt{q}]$, $d=m^2-4q$.
Let $A$ be an abelian variety over $\F_q$ isogenous to
$E_m \times E_m \times E_{m-2}$ which has 
an indecomposable principal polarization.  Then there exists
a rank $2$ indecomposable positive definite hermitian form
of discriminant $2$ on $\Z[\pi]$, $\pi = \frac{-m + \sqrt{d}}{2}$.
\end{thm}
\noindent
{\bf Proof:} 
We have that $m=2x+1$ and Propositions \ref{primeto} and \ref{primetom-2} 
show that the assumptions of the theorem imply that $m$ and $m-2$ are
prime to the characteristic. In fact, the restriction $p \ne 5$
is not necessary if there are no solutions to the equation 
$5^e=x^2+x+3$ satisfying $e \ge 3$, $e$ odd and $5$ dividing $2x-1$.

So all abelian varieties in this proof are ordinary.  Due to
Fact 4.1, we also know that $2$ is prime to the characteristic,
so the group schemes we work with will be \'etale.
Thus we will identify finite group schemes with their 
$\overline{\F_q}$-points.
In the cases at hand, we have $d= -3$ or $d=-11$, so in particular,
$d$ is the discriminant of an imaginary quadratic field of class number one.

By the equivalence of categories given in the Appendix, it
suffices to show that the indecomposable principal polarization
on $A$ induces an indecomposable polarization on $E_m \times E_m$
of degree $2$, where $E_m$ is an elliptic curve with $q+1-m$ points
over $\F_q$.

Let $F$ be the Frobenius endomorphism on $A$ and $V$ the Verschiebung, 
so that $\phi = F+V$ has eigenvalues $-m,-m,-m,-m,-(m-2),-(m-2)$. 
Let $$B= \text{ the connected component of the kernel of  } \phi+m,$$  and  
$$C=\text{ the connected component of the kernel of  } \phi+m-2.$$  
Then since $\Q(\sqrt{d})$ has class number one,
$$B \simeq E_m \times E_m$$ and $$C \simeq E_{m-2}.$$ 
We have $A$ is isogenous to $B \times C$ 
$$f:(B \times C)  \rightarrow A,$$ 
with kernel $\Delta$ isomorphic to the intersection of $B$ and $C$, 
$\Delta \simeq  B \cap C$. The polarization 
$\lambda$ on $A$ induces polarizations $b$ on $B$ and $c$ on $C$.
It suffices to show that $b$ is indecomposable of degree $2$. 

Let $N_b ={\rm Ker}(b)$,  $N_c ={\rm Ker}(c)$.
Together $(b,c)$ is the induced polarization on $B \times C$,
with kernel $N_b \times N_c$, and so 
we must have $$ B \cap C \subset N_b \times N_c,$$
embedded diagonally.  We want to show that $N_b$ has order $4$.

By the definition of $B$ and $C$, we know that
$B \cap C$ is killed by $2$, so it is contained in $C[2]$,
which has order $4$.  In addition, $B \cap C$ is stable
under the action of Frobenius, so it is an $R$-module,
where $R = \Z[\pi] = \Z[x]/(x^2-mx+q)$.  But $(2)$ is
inert in $R$, so $R/2R \simeq \F_4$.  So $B \cap C$
is an $\F_4$-vector space of order less than or equal to $4$.
$B \cap C$ cannot be trivial, or else $A$ would be split and
the polarization would be decomposable. So $B \cap C$ must
have order $4$ and the polarization $(b,c)$ has degree $4$.
Since $B \cap C$ is embedded diagonally into $N_b \times N_c$,
it follows that $N_b$ has order $4$.

In addition, $B$ is orthogonal to $C$ with respect to the
skew-symmetric bilinear pairing $\langle, \rangle$
defined by the polarization of $A$ on the $\ell$-adic
Tate module associated to $A$.  
This follows from the fact that $\phi$ is hermitian with respect to the
pairing, since 
$$\phi^* = (F+V)^*=F^*+V^*=V+F=\phi.$$
Thus for elements $b \in T_lB$ and $c \in T_lC$ of the $\ell$-adic
Tate modules associated to $B$ and $C$, we have
$$ 2\langle b, c \rangle = \langle b, 2c \rangle = 
\langle b, (\phi + m)c \rangle = \langle (\phi + m)b, c \rangle = 0.$$
This shows that they are orthogonal since $\Z_{\ell}(1)$ is torsion-free.
Thus $B \cap C$ is maximal isotropic with respect to the pairing.   

Finally, the polarization $b$ is indecomposable, else $\lambda$ would not be.
$\square$

\begin{cor} There are no defect $2$ curves in genus $3$ over $\F_q$
if $q$ is of the form $q=x^2+x+a$, $a=1$ or $a=3$ with $a \le x$.
\end{cor}
{\bf Proof:}  This is a direct consequence of Propositions
\ref{defect2} and \ref{HM2} and Theorem \ref{notglue2}.
The exceptional cases can be handled individually. If $q=3$,
there is no defect $2$ curve (see Section 4.1.1).
If $q=343$, then $m-2=35$, so there is no elliptic curve over
$\F_q$ with trace $\pm (m-2)$. Thus there is no curve of
type $[m,m,m-2]$ over $\F_{343}$.  Similiarly, if a solution to the equation 
$5^e=x^2+x+3$ exists with $e \ge 3$, $e$ odd and 
$5$ dividing $2x-1$, then a curve of type $[m,m,m-2]$ over
$\F_{5^e}$ would fail to exist for the same reason.

\vspace{10pt}
Thus we have established that $N_q(3) \le q+1+3m-3$ for such $q$.
To finish the proof of Theorem \ref{d-3-11}, we must show that 
there exists a curve $C$ of genus $g(C)=3$ over $\F_q$ such that 
$|N(C)-(q+1)|=3m-3.$

\begin{prop} \label{m-1}
Let $q=p^e=x^2+x+a$, $a=1$ or $a=3$, with $a \le x$, and
$q \ne 3$, $q \ne 3^5$.  Then there exists a curve $C$ of genus $g(C)=3$ 
over $\F_q$ with zeta function of type $\pm [m-1,m-1,m-1]$.
\end{prop}
{\bf Proof:}  Define $d'=(m-1)^2-4q$.  Recall that 
$$d=m^2-4q=-3 \text{ or } -11.$$
We first show that $(m-1)$ is prime to $p$ and 
$$d' \notin \{-3,-4,-8,-11\}.$$

To check that $(m-1)$ is prime to $p$,
write $m-1=2x$. By Fact 4.1 we know that $p \ne 2$.  So if 
$p$ divides $(m-1)$ then $p$ must divide $x$.  
For $p^e=x^2+x+1$ this is impossible. 
If $p^e=x^2+x+3$ and $p$ divides $x$, then $p=3$. 
As explained in Proposition \ref{primetom-2} above,
the only powers of $3$, $e$ odd, of the
form  $3^e = x^2 + x + 3$  are  $3$  and  $3^5 = 243$. 

To show that $d' \notin \{-3,-4,-8,-11\}$
we write $d'=d-2m+1$, where $d=-3$ or  $d=-11$.  If $d=-3$, then
$$d'=-2m-2.$$ So $d' \in \{-3,-4,-8,-11\}$ only if $m=1$ or $m=3$, which
occurs only for $q=3$, $m=3$.  Indeed, we assumed that $q\ne 3$, 
since in that case $[m-1,m-1,m-1]$ is
not possible.  If $d=-11$, then $$d'=-2m-10,$$ which is not in the set
$\{-3,-4,-8,-11\}$ for any $m>1$. 

Now we can apply the theory of Hermitian modules.
By Theorem 8.2 in \cite{Hoffmann}, there exists an 
indecomposable, positive definite, unimodular Hermitian module
over $R_{d'}$  of rank $3$.  Applying the functor $S$ from the Appendix
to this module, we obtain an abelian  variety $A$ isogenous to
$E_{m-1} \times E_{m-1} \times E_{m-1}$ with an indecomposable principal
polarization.  Then by the Torelli theorem (cf. \cite{LauterJAG}, Appendix), 
there exists a genus $3$ curve $X$ over $\F_q$ whose Jacobian is isomorphic
either to $A$ or to the quadratic twist of $A$; hence $X$ is of type 
$\pm [m-1,m-1,m-1]$.  $\square$

\vspace{10pt}
\noindent {\bf Note 4.1.4}
Indeed if $q=243$, there is no curve of type $\pm [m-1,m-1,m-1]$
since $m-1=30$, and Honda-Tate theory shows
that there is no abelian variety of dimension $3$ over $\F_{243}$ corresponding
to that trace of Frobenius.

\begin{cor} For $q=7$ or $q=13$, we have $N_q(3)=q+1+3m-3$.
\end{cor}
{\bf Proof:}  This was  noted in \cite{Serre-Harvard}.  It follows
from Proposition \ref{m-1} and the fact that a Frobenius with the opposite
sign is not possible because the corresponding curve would have a 
negative number of points.


\begin{prop} Let $q=3^5=243$.  Then there exists a curve $C$ of genus $3$ over
$\F_{q}$ with $|N(C)-(q+1)|=3m-3$.
\end{prop}
Proof: The only possible zeta function for defect $3$ 
in this case is $[m,m,m-3]$.  Since $m=31$ is prime to $3$, 
there exist elliptic curves $E_m$ and $E_{m-3}$. Take the
polarization $b$ on $B=E_m \times E_m$ of discriminant $3$
given by the matrix
\begin{displaymath}
       \left[ \begin{array}{cc}
           2 & 1\\
           1 & 2 \\
      \end{array} \right].
    \end{displaymath}
For the polarization on $C=E_{m-3}$ we take $3$ times the canonical 
polarization.  To glue $B$ to $C$ we must find an isomorphism of
the kernels of $b$ and $c$ {\it as group schemes}, and since the 
order is not prime to the characteristic, it is not enough to consider
$\overline{\F_q}$-points.  The kernel of $B$ is isomorphic to the
$3$-torsion of $E_m$.  So both kernels are of (\'etale, local) type:
$\Z/3\Z \times \mu_3$(twisted), and there is only one choice for the 
pairing up to sign.
$\square$

\vspace{10pt}
\noindent
{\bf Note 4.1.5} An explicit search for a genus $3$ curve
over $\F_{243}$ with the maximum or minimum number of points currently
seems out of reach.  Even counting the points on each of the
approximately $243^6$ homogeneous plane quartics in the variables
$x^2$, $y^2$, $z^2$ using the naive method would take $3^{40}$ steps.
A computation of this size is currently infeasible.
A more fruitful approach may be to generate the equations of the elliptic
curves $E_m$ and $E_{m-3}$, which we can easily do, and to try to 
glue them together explicitly.
  
\subsection{Proof of Theorem \ref{d-4-8}}

Suppose $q=x^2+b$, $b=1$ or $b=2$ with $b$ 
satisfying $b \le x$.
Again by Fact 3.1 there are no defect $0$ curves since 
$d=-4$ or $-8$.  By Fact 3.2, defect $1$ curves do not exist for
$g>2$.  It follows that $$N_q(3) \le q+1+3m-2.$$ 
It remains to show that there exists
a curve $C$ over $\F_q$ such that $|N(C)-(q+1)|=3m-2$.

\begin{prop} \label{m-2}
Suppose $q=x^2+b$, $b=1$ or $b=2$ with $b$ 
satisfying $b \le x$.
Then $(m,p)=1$ and $(m-2,p)=1$ unless $q=2$.
\end{prop}
{\bf Proof:}  By Fact 4.1, $p$ is not equal to $2$ unless $e=1$.  
Suppose that $q \ne 2$.  
Write $$m=2x \quad \text{ and  } \quad m-2=2(x-1).$$  
If $p$ divides $m$, then since $p \ne 2$ we
have $p$ divides $x$, which is not possible for  either $b=1$ or $b=2$.  

If $p$ divides $(m-2)$, then $p$ must divide $(x-1)$.  
If $b=1$, then this is impossible since
$$p^e=(x-1)(x+1)+2$$ and $p \ne 2$. If $b=2$ and $p$ divides $(x-1)$, 
then $p$ divides $(q-3)$, so
$p=3$.  But there are no solutions to the equation 
$$q=3^e=x^2+2,$$
meeting the congruence condition $$x \equiv 1 \pmod{3}$$ with $e>1$, since by
\cite{Apery} there exist at most two solutions to such an equation
and in this case they are $e=1$, $x=1$ and $e=3$, $x=5$.
The second solution does not satisfy the congruence condition
and the first solution does not matter since we do not write
$q=3$ in the form $x^2+2$ (See Note 2.2).

\subsubsection{q=2}
When $q=2$, a curve $C$ of genus $3$ and defect $2$ over $\F_2$ 
was given in \cite{Serre-Harvard}.  Its equation is:
$$ x^3y+y^3z+z^3x + x^2y^2+y^2z^2+z^2x^2+x^2yz+y^2xz = 0$$
and it has $7$ points, which is defect $2$ since $m=2$.
It has zeta function of type 
$$[m+1-4\cos^2(\frac{\pi}{7}), m+1-4\cos^2(\frac{2\pi}{7}), m+1-4\cos^2
(\frac{3\pi}{7})].$$
This is possible because $$2\sqrt{2} > 1-4\cos^2(\frac{3\pi}{7}),$$
and it is the only zeta function possible for this case.
It is a twist of the Klein curve which becomes isomorphic to the Klein curve
over the field $\F_{2^7}$.

It is interesting to note that when $q=2$, the zeta function
type $[m, m, m-2]$ is not possible because the curve would have $7$
points over $\F_2$ but only $1$ point over $\F_8$.  Alternatively,
we can show that the glueing of the supersingular elliptic curves
is not possible by examining the group schemes in question and 
showing there is no isomorphism between them.

The zeta function type $[m, m+\sqrt{3}-1, m-\sqrt{3}-1]$ 
is not possible because
the curve would have  $7$ points over $\F_2$ but only $5$ points over $\F_4$.

\subsubsection{$q \ne 2$}
Now assume $q \ne 2$.  In order to have a
genus $3$ curve of type $[m,m,m-2]$, there must exist an indecomposable
Hermitian module of rank $2$ and discriminant $2$ over $R_d$.
When $d=-4$ or $d=-8$, an indecomposable Hermitian module of rank $2$ and
discriminant $2$ exists and is given for $d=-4$ by 
\begin{displaymath}
       \left[ \begin{array}{cc}
           2 & 1+i \\
           1-i & 2 \\
      \end{array} \right]
    \end{displaymath}  
and for $d=-8$ by
\begin{displaymath}
       \left[ \begin{array}{cc}
           2 & -\sqrt{-2}\\
        \sqrt{-2} & 2 \\
      \end{array} \right].
    \end{displaymath}  
In both cases, the module is indecomposable because all 
the values taken by the Hermitian form are divisible by $2$, so $1$ is not
represented.  But if the form were equivalent to 
\begin{displaymath}
       \left[ \begin{array}{cc}
           1 & 0 \\
           0 & 2 \\
      \end{array} \right],
    \end{displaymath}
then $1$ would be represented.

The next theorem shows that in this case it is possible to glue two abelian
varieties together to obtain the Jacobian of a genus $3$ curve with the
desired property.

\begin{thm} \label{glue2} 
Suppose $q=x^2+j$, $j=1$ or $j=2$ with $j$ 
satisfying $j \le x$, and suppose $q \ne 2$.
Then there exists an abelian variety $A$ over $\F_q$,
isogenous to $E_m \times E_m \times E_{m-2}$,
with an indecomposable principal polarization.
\end{thm}
{\bf Proof:}  By Fact 4.1, the characteristic of the field is not equal
to $2$, and by Proposition \ref{m-2}, there exist abelian 
varieties $B=E_m \times E_m$ and $C=E_{m-2}$.

{\bf j=1.} Let $b$ be the polarization on $B$ corresponding to the 
positive definite indecomposable Hermitian module 
\begin{displaymath}
       \left[ \begin{array}{cc}
           2 & 1+i \\
           1-i & 2 \\
      \end{array} \right].
    \end{displaymath}
Let $c$ be two times the canonical polarization on $C$.
It has kernel equal to the $2$-torsion of $C$.
We proceed by calculating the kernel of $b$ and then glueing
it to the kernel of $c$.  The order of (the $\overline{\F_q}$ points of) 
these group schemes is $4$, which is prime to the characteristic. 
By Mumford's criteria, we need to find an isomorphism of the Galois modules
which is an anti-isometry with respect to the pairings.
If necessary, we will replace $C$ by an isogenous elliptic curve.

The kernel of $b$ is $E_m[\lambda]\times E_m[\lambda]$, 
where $E_m[\lambda]$ is the $\lambda$-torsion of $E_m$,
and $\lambda = 1+i$.  The $\lambda$-torsion of the elliptic curve
is contained in the $2$-torsion, since $2=(1+i)(1-i)$.
The Frobenius of $E_m$ acts on the $2$-torsion of $E_m$ by fixing
one of the three non-trivial points and by exchanging the other two.
This can be seen by looking at the  $2 \times 2$ matrix which represents
the action of Frobenius on the $\ell$-adic Tate module when $\ell=2$.
The characteristic polynomial of Frobenius is $t^2+mt+q,$ where $-m$
is the trace of the matrix and $q$ is the determinant.  In this case,
$q=x^2+1$ and $m=2x$, with $x$ even since $q$ is odd.  Thus we have
$m \equiv 0 \pmod{4}$.  So the number of points on the elliptic curve
over $\F_q$, $N(E_m)$, satisfies
$$N(E_m) = q+1+m \equiv 2 \pmod{4}.$$
Similarly,
$$N(E_{m-2}) = q+1+m-2 \equiv 0 \pmod{4}.$$
In both cases, the trace is even and the determinant is odd, 
so the possibilities for the matrix of Frobenius $\pmod{2}$ are:
\begin{displaymath}
       \left[ \begin{array}{cc}
           1 & 0 \\
           1 & 1 \\
      \end{array} \right],
   \left[ \begin{array}{cc}
           1 & 1 \\
           0 & 1 \\
      \end{array} \right],
  \left[ \begin{array}{cc}
           0 & 1 \\
           1 & 0 \\
      \end{array} \right], \text{ and }
  \left[ \begin{array}{cc}
           1 & 0 \\
           0 & 1 \\
      \end{array} \right].
    \end{displaymath}
The first three of these matrices act on the $2$-torsion by fixing
one of the three non-trivial points and by exchanging the other two.
The last one is the identity matrix and fixes all three non-trivial 
$2$-torsion points.
But the identity matrix is only possible if $N(E) \equiv 0 \pmod{4}$,
while the first three are possible in either case.
This follows by writing the matrix with entries
\begin{displaymath}
       \left[ \begin{array}{cc}
           a & b \\
           c & d \\
      \end{array} \right],
      \end{displaymath}
the number of points on $E$ as $$N(E)=(ad-bc)+1-(a+d),$$
and considering all the possibilities for the entries $\pmod{4}$.

So we only need to chose $C$ to be an elliptic curve $E_{m-2}$ with
Frobenius fixing one of its $2$-torsion points and switching the other
two.  It suffices to take 
an elliptic curve with endomorphism ring $\Z[\pi] \subset R'$
such that $$\frac{\pi -1}{2} \not\in R'.$$  
It follows from Deuring that such a curve exists.
By chosing $C$ in this
way we can find an isomorphism of $C[2]$ with the kernel
of $b$ which respects the action of Frobenius.  Finally,
the Galois module isomorphism must be an anti-isometry with
respect to the pairings.  The pairings take values $\pm 1$.

{\bf j=2} The proof is almost the same except in this case 
$\lambda = \sqrt{2}$ and the polarization on $B$ is given by
the matrix
\begin{displaymath}
       \left[ \begin{array}{cc}
           2 & -\sqrt{-2}\\
        \sqrt{-2} & 2 \\
      \end{array} \right].
    \end{displaymath}  
Again we have 
$$N(E_m) = q+1+m = x^2+3+2x \equiv 2 \pmod{4},$$
and
$$N(E_{m-2}) = q+1+m-2 =x^2+1+2x \equiv 0 \pmod{4},$$
since $x$ is odd. So the same argument works. $\square$

\begin{cor}
If $q=x^2+j$, $j=1$ or $j=2$ with $j$ 
satisfying $j \le x$, and $q \ne 2$, then
there exists a curve of type $\pm [m,m,m-2]$.
\end{cor}

Proof: This follows from Theorem \ref{glue2} and the fact that 
(see \cite{Oort}): if  $A$  is a principally polarized indecomposable 
abelian variety over an algebraically closed field, of dimension $3$, 
then $A$  is the Jacobian of a curve. The ``precise Torelli''
theorem in the Appendix to \cite{LauterJAG} allows you to descend 
from the algebraically closed field to any field at the cost of
quadratic twist, as is explained there.

\vspace{10pt}
\noindent
{\bf Example} An example of a genus $3$ curve corresponding to this type of
glueing was found by van der Geer and van der Vlugt.  It can be described as
the curve over $\F_{27}$ obtained by taking the fiber product of
$$y^2=x^3+2x^2+2x$$ with
$$y^2=2x^3+2\alpha^4x^2+\alpha^8x$$
where $\alpha^3+2\alpha^2+1=0$.    
It is a defect $2$ curve with $56$ points (instead of $58$ as stated
in \cite{Geer}).

\subsection{Proof of Theorem \ref{all}}

Let $q=p^e$ be a power of a prime.
We divide the proof into two cases: $e$ even and $e$ odd.

{\bf e even.} First suppose $e=2r$ is even.  
Then Ibukiyama has shown \cite{Ibu}
that, for $p$ an odd prime, there exists a curve $C$ of genus
$3$ over $\F_p$ such that 
$$\#C(\F_{p^{2r}}) = 1 + p^{2r} + (-1)^{r+1}6p^r.$$
Thus for all $q$ an even power of an odd prime, there exists
a genus $3$ curve attaining either the Weil maximum or the 
Weil minimum.  

Now suppose that $p=2$, so that $q=2^{2r}$.  Then $m=2^{r+1}$,
and $(m-1)$ is prime to $p$.  If we let $d'=(m-1)^2-4q$, 
then $$d'=-2^{r+2}+1,$$
and so we see that for $r \ge 1$, $$d' \notin \{-3,-4,-8,-11\}.$$
By Theorem 8.2 in \cite{Hoffmann}, there exists an 
indecomposable, positive definite, unimodular Hermitian module
over $R_{d'}$  of rank $3$.  Applying the functor $S$ from the Appendix
to this module, we obtain an abelian  variety isogenous to
$E_{m-1} \times E_{m-1} \times E_{m-1}$ with an indecomposable principal
polarization.  Then by the Torelli theorem, (see the Appendix 
to \cite{LauterJAG}) there exists a genus $3$
curve over $\F_q$ which is of type $\pm [m-1,m-1,m-1]$.  

For example, over $\F_4$,  the equation of
the Klein curve is given in \cite{Serre-Harvard} to show that
$N_4(3)=14$.  It is a curve of type $[m-1,m-1,m-1]$.  

{\bf e odd.} Now suppose that $e$ is odd.  As usual, write $$q=x^2+x+a,$$ with
$-x \le a \le x$.  We divide the proof according to the
value of $d=m^2-4q$.  If $d \in \{-3,-4,-8,-11\}$, then these cases have been
treated in Theorems \ref{d-3-11} and \ref{d-4-8}.  

If $d$ is not in the set$\{-3,-4,-8,-11\}$,  we check whether $m$
is prime to $p$.  If it is, then a curve of type $[m,m,m]$ or 
$[-m,-m,-m]$ exists.  If not, then $(m-1)$ is prime to $p$.  
It remains to check that 
$$d'=(m-1)^2-4q \notin \{-3,-4,-8,-11\}.$$ 
Since $$d'=d-2m+1$$ and $$d \equiv 0 \text{ or  }1 \pmod{4},$$ this could
only occur for $d=-7$ and $m=1$, but $m>1$, so it does not occur. Thus there
exists a curve of type $[m-1,m-1,m-1]$ or of type $[-(m-1),-(m-1),-(m-1)]$ in
this case.

\clearpage

{\large \bf{Appendice \quad \quad \quad \quad  J-P.Serre, 1999}}

\begin{center}
{\bf Modules hermitiens et courbes alg\'ebriques}
\end{center}

\bigskip
\noindent {\bf 1. Notations}

On note $k$ un corps fini \`a $q$ \'el\'ements de caract\'eristique $p$.  On se
 donne un entier $a$ et l'on pose $d=a^2-4q$.  On suppose:
\begin{description}
\item[(1.1)] {\it $a$ est premier \`a $p$;}
\item[(1.2)] {\it $d<0$;}
\item[(1.3)] {\it $d$ est le discriminant d'un corps quadratique imaginaire.}
\end{description}

On pose $R=\Z[X]/(X^2-aX+q)$.  Vu les hypoth\`eses ci-dessus, $R$ est
l'anneau des entiers du corps quadratique imaginaire $K=\Q(\sqrt{d})$,
dans lequel $p$ est d\'ecompos\'e.

(L'hypoth\`ese (1.3) n'est past indispensable pour la suite; il est
souvent commode de ne pas la faire; on doit alors travailler avec des
``ordres'' non maximaux de $K$.)

On choisit une courbe elliptique $E$ sur $k$ dont les valeurs propres
de Frobenius sont $\frac{1}{2}(a\pm\sqrt{d})$ (son nombre de points
est donc $q+1-a$).  On sait qu'il en existe.

\bigskip
\noindent
{\bf 2. La cat\'egorie $\Ab(a,q)$}

On note $\Ab(a,q)$ la cat\'egorie des vari\'et\'es ab\'eliennes sur $k$
ayant les propri\'et\'es \'equivalentes suivantes:
\begin{description}
\item[(2.1)] {\it $A$ est $k$-isog\`ene \`a un produit de copies de $E$.}
\item[(2.2)] {\it Si $F_A$ et $V_A$ d\'esignent respectivement le Frobenius
et le ``Verschiebung'' de $A$, on a $F_A+V_A=a$ dans $\End(A)$.}
\item[(2.3)] {\it Les valeurs propres de $F_A$ sont celles de $F_E$,
r\'ep\'et\'ees $g$ fois o\`u $g=\dim A$.}
\end{description}
(L'\'equivalence de ces propri\'et\'es r\'esulte de th\'eor\`emes de
Tate.)  On aura besoin plus loin de la propri\'et\'e suivante:
\begin{description}
\item[(2.4)] {\it Si $A,B$ appartiennent \`a $\Ab(a,q)$, et si
$f:A\rightarrow B$ est un homomorphisme d\'efini sur une extension
$k'$ de $k$, alors $f$ est d\'efini sur $k$.}
\end{description}
En effet, on peut supposer que $k'$ est une extension finie de $k$.
Soit $m\ge 1$ son degr\'e.  Puisque $f$ est d\'efini sur $k'$, il
commute \`a la puissance $m$-\`eme du Frobenius.  Mais si $\pi$ est le
Frobenius (i.e. le g\'en\'erateur ``$X$'' de $R$), l'anneau $\Z[\pi^m]$
est un sous-anneau {\it d'indice fini} de $R=\Z[\pi]$.  La commutation avec
$\pi^m$ entra\^{i}ne donc celle avec $\pi$.

\bigskip
\noindent
{\bf 3. Une \'equivalence de cat\'egories}

Notons $\Mod(R)$ la cat\'egorie des $R$-modules sans torsion de type
fini (i.e. projectifs de type fini, vu que $R$ est un anneau de
Dedekind, gr\^{a}ce \`a (1.3)).

Noter que $R$ op\`ere sur toute vari\'et\'e $A$ de $\Ab(A,q)$; de
plus, d'apr\`es (2.4), tout $f:A\rightarrow B$ est un
$R$-homomorphisme.

Si $A\in \Ab(a,q)$, notons $T(A)$ le $R$-module $\Hom(E,A)$.  C'est un
\'el\'ement de $\Mod(R)$.

\begin{description}
\item[(3.1)] {\it Le foncteur $T:\Ab(A,q)\rightarrow\Mod(R)$ est une
\'equivalence de cat\'egories.}
\end{description}
On a en particulier un isomorphisme naturel:
$$ \Hom(A,B) = \Hom_R(T(A),T(B)).$$
On peut expliciter un foncteur $S:\Mod(R)\rightarrow\Ab(a,q)$ qui est
``inverse'' au foncteur $T$: si $L\in\Mod(R)$, on d\'efinit $S(L)$
comme la vari\'et\'e ab\'elienne $L\otimes_R E$, ``produit tensoriel''
de $L$ par $E$ (de tels produits tensoriels existent dans toute
cat\'egorie ab\'elienne; on \'ecrit $L$ comme conoyau d'un
homomorphisme $R^N\rightarrow R^M$ et l'on d\'efinit $L\otimes_R E$
comme le conoyau de l'homomorphisme correspondant: $E^N\rightarrow
E^M$).

Les assertions ci-dessus entra\^{i}nent en particulier:
\begin{description}
\item[(3.2)] {\it Toute vari\'et\'e ab\'elienne $A$ appartenant \`a $A(a,q)$
peut s'\'ecrire sous la forme $L\otimes_R E$ avec $L=\Hom(E,A)$.}
\end{description}
(Dans ce qui suit, j'abr\`egerai $L\otimes_R E$ en $A_L$.)

Bien s\^{u}r, on a:
\begin{description}
\item[(3.3)] $\rang(L) = \dim A_L$.
\item[(3.4)] $A_R = E$.
\end{description}
Si $M$ est un sous-module d'indice fini de $L$, l'inclusion
$M \rightarrow L$ d\'efinit un morphisme $f:A_M \rightarrow A_L$ qui est
une {\it isog\'enie}; de plus:
\begin{description}
\item[(3.5)] {\it Le degr\'e de $f$ est \'egal \`a $(L:M)$.}
\end{description}
\textbf{Exemple:} {\it le cas o\`u $\dim A=1$.}

Ce cas correspond \`a $\rang(L) = 1$; autrement dit $L$ est un
$R$-module inversible.  Les classes de tels modules correspondent aux
\'el\'ements de $\Cl(R)=\Pic(R)$; leur nombre est le {\it nombre de classes}
$h(d)$ de l'anneau $R$.  On conclut de l\`a que le nombre des classes
d'isomorphisme de courbes elliptiques $k$-isog\`enes \`a $E$ est
\'egal \`a $h(d)$; on retrouve un r\'esultat bien connu.

\bigskip
\noindent{\bf 4. Dualit\'e}

Si $A$ appartient \`a $\Ab(a,q)$ il en est de m\^{e}me de sa {\it duale}
$A^*$, puisque $A$ et $A^*$ sont $k$-isog\`enes.

Si $L$ appartient \`a $\Mod(R)$, notons $R^*$ son {\it anti-dual}, autrement
dit l'ensemble des homomorphismes $f:L\rightarrow R$ qui sont
anti-lin\'eaires, i.e. tels que $f(rx) = \overline{r}f(x)$ pour $r\in
R$ et $x\in L$.

Les foncteurs $A \mapsto A^*$ et $L \mapsto L^*$ se correspondent par le
dictionnaire du $\S$3.  Autrement dit:
\begin{description}
\item[(4.1)] {\it Si $A\in \Ab(a,q)$ correspond au module $L$, sa duale
$A^*$ correspond au module $L^*$.}
\end{description}
(Le fait que l'on prenne l'{\it anti}-dual, et non le dual, provient de ce
que la transposition sur $\End(E)=R$ est la conjugaison complexe.)

\bigskip
\noindent
{\bf 5. Polarisations}

Une polarisation est un morphisme $\varphi: A\rightarrow A^*$ qui
provient d'un diviseur ample sur $A$ (cf. [14]).  Si
$A$ correspond au module $L$, $\varphi$ correspond par (4.1) \`{a} un
morphisme $L\rightarrow L^*$, autrement dit \`a une forme
sesquilin\'eaire $H:L\times L\rightarrow R$.  De plus:
\begin{description}
\item[(5.1)] {\it $H$ est une forme hermitienne d\'efinie $>0$.}
\end{description}
(Autrement dit on a $H(x,y) = \mbox{conjugu\'e de\ }H(y,x)$ pour
$x,y\in L$, et $H(x,x)>0$ si $x\not=0$.)

Inversement:
\begin{description}
\item[(5.2)] {\it Toute forme hermitienne d\'efinie $>0$ sur $L$ d\'efinit
une polarisation de $A$.}
\end{description}
Une polarisation $\varphi$ a un degr\'e $\deg(\varphi)$ d\'efini par:
$$ \deg(\varphi)^2 = \mbox{ordre du sch\'ema en groupes fini
$\Ker(\varphi)$}.$$
En termes de $L$ et de la forme hermitienne $h:L\rightarrow L^*$, ceci
se traduit par:
\begin{description}
\item[(5.3)] $\deg(\varphi)^2 = (L^*:hL)$.
\end{description}
En particulier:
\begin{description}
\item[(5.4)] {\it Pour que $\varphi$ soit une polarisation principale
(i.e. $\deg(\varphi)=1$), il faut et il suffit que $hL=L^*$} (auquel
cas on dit que $H$ est $R$-non d\'eg\'en\'er\'ee, ou encore que le
discriminant du module hermitien $L$ est \'egal \`a 1.)
\end{description}

Lorsque $L$ est un $R$-module libre (ce qui est toujours le cas si
$h(d)=1$) et qu'on en choisit une base $(e_i)$, la forme $H$ est
donn\'ee par une matrice hermitienne $(r_{ij})$ et la condition que le
discriminant soit $1$ se traduit par $\det((r_{ij}))=1$.

\bigskip
\noindent
{\bf 6. Polarisations principales ind\'ecomposables}

Soit $L\in\Mod(R)$, muni d'une forme hermitienne $>0$ de discriminant
$1$.  Soit $A$ la vari\'et\'e ab\'elienne polaris\'ee correspondante.
Vu le $\S 5$, on a:
\begin{description}
\item[(6.1)] {\it Pour que $A$ soit ind\'ecomposable} (comme vari\'et\'e
ab\'elienne polaris\'ee), {\it il faut et il suffit que $L$ soit 
ind\'ecomposable comme module hermitien.}
\end{description}

Noter que, \`a cause de (2.4), la notion d'ind\'ecomposabilit\'e pour
$A$ a le m\^{e}me sens sur $k$, ou sur toute extension de $k$.  Il
s'ensuit qu'un module hermitien ind\'ecomposable donne une vari\'et\'e
ab\'elienne polaris\'ee qui est ``absolument'' ind\'ecomposable.

On peut donner des exemples d'anneaux $R$ qui n'ont aucun module
hermitien ind\'ecomposable (de discriminant 1) en dimension $g$
\'egale \`a 2 ou 3.  D'apr\`es Hoffmann (cf [5]), ce sont:
\begin{description}
\item[(6.2)] {\it Pour $g=2$, les cas $d=-3,-4,-7$;}
\item[(6.3)] {\it Pour $g=3$, les cas $d=-3,-4,-8,-11$.}
\end{description}
Hoffmann a \'egalement montr\'e que, pour $g=2,3$ les valeurs ci-dessus
sont les {\it seules valeurs} de $d$ pour lesquelles tous les modules
hermitiens de discriminant $1$ sont d\'ecomposables.

\bigskip
\noindent
{\bf 7. Application aux courbes alg\'ebriques}

Ce qui pr\'ec\`ede s'applique \`a la jacobienne $J$ d'une courbe de
genre $g$ d\'efinie sur $k$ (et dont les valeurs propres de Frobenius
sont celles de $E$ r\'ep\'et\'ees $g$ fois).  Comme $J$ est munie
d'une polarisation principale ind\'ecomposable, on d\'eduit de l\`a:
\begin{description}
\item[(7.1)] {\it Pour $g=2,3$ il n'existe aucune courbe $C$ dont la
jacobienne appartienne \`a $\Ab(a,q)$, si $d=a^2-4q$ a l'une des
valeurs donn\'ees dans (6.2) et (6.3).}
\end{description}

Posons, comme d'habitude $m=[2q^{1/2}]$ et supposons $m\not=0 \pmod
p$.  Supposons que le nombre $N$ de points de la courbe soit \'egal
\`a $q+1+gm$.  On sait que $J$ appartient alors \`a $\Ab(a,q)$, avec
$a=-m$.  On d\'eduit de l\`a et de (7.1) que la courbe en question
n'existe pas lorsque:

$g = 2$, $\quad$ $m^2-4q=-3, -4, \mbox{ou\ }-7$;

$g = 3$, $\quad$ $m^2-4q=-3,-4,-8, \mbox{ou\ }-11.$

(M\^{e}me r\'esultat si $N=q+1-gm$, car cela ne fait que changer le
signe de $a$.)

Ainsi, par exemple, il n'existe aucune courbe de genre $3$ sur
$\mbox{\textbf{F}}_{27}$ ayant 58 points, car une telle courbe
donnerait $d=-8$.

On peut aussi proc\'eder en sens inverse, et utiliser des modules
hermitiens ind\'ecomposables de rang 2 ou 3 pour construire (ou
plut\^{o}t pour prouver l'existence\dots) de courbes.  En effet, soit
$L$ un $R$-module hermitien ind\'ecomposable $R$-non d\'eg\'en\'er\'e
de rang 2 (resp. 3).  En appliquant le th. de Torelli \`a $A_L$, on en
d\'eduit:
\begin{description}
\item[(7.2)] {\it Si $g=2$, il existe une courbe $C$ sur $k$ dont la
jacobienne est isomorphe \`a $A_L$} (et qui a donc $q+1-2a$ points).
\item[(7.3)] {\it Si $g=3$, il existe une courbe $C$ sur $k$ dont la
jacobienne est isomorphe, soit \`a $A_L$, soit \`a la ``tordue
quadratique'' de $A_L$} (et qui a donc $q+1-3a$ points dans le premier cas
et $q+1+3a$ points dans le second cas).
\end{description}

De plus, dans le cas $g=3$, si la courbe consid\'er\'ee n'est pas
hyperelliptique le cas ``$-3a$'' exclut le cas ``$+3a$''.  Autre propri\'et\'e
de ce cas (cons\'equence de Torelli, ici aussi): si $C$ est la courbe
consid\'er\'ee, suppos\'ee non hyperelliptique, son groupe
d'automorphismes $\Aut(C)$ est un sous-groupe d'indice 2 de $\Aut(L)$; de 
fa\c con plus pr\'ecise, on a
$$ \Aut(L) = \{\pm 1\} \times \Aut (C).$$
Par contre, si $g=2$, ou si $g=3$ et $C$ est hyperelliptique, on a
$$ \Aut(L) = \Aut(C).$$

{\bf Exemple:} Prenons $q=41$, $g=3$, $a=\pm m = \pm 12$, de sorte que
$d=-20$.  D'apr\`es Hoffmann, il y a deux possibilit\'es pour $L$,
avec chaque fois $\Aut(L)=\{\pm 1\}\times S_3$.  D'o\`u l'existence de
courbes de genre 3 sur $\mbox{\bf F}_{41}$, ayant soit 78 points (ce qui
serait le maximum), soit 6 points (ce qui serait le minimum).  Chacune
de ces deux courbes a un groupe d'automorphismes qui est, soit $\{\pm
1\}\times S_3$, soit $S_3$.

{\bf Remarque:} On peut utiliser (7.2) et (7.3) pour d\'emontrer
certains cas de (6.2) et (6.3).  Prouvons par exemple que, pour
$d=-8$, $g=3$, il n'y a pas de module hermitien ind\'ecomposable de
discriminant 1.  S'il y en avait un, par (7.3) appliqu\'e \`a $q=3$,
$a=2$, il y aurait: 

soit une courbe $C$ dont les valeurs propres de
Frobenius sont $1\pm\sqrt{-2}$ (r\'ep\'et\'ees 3 fois): son nombre de
points serait $q+1-3a=3+1-6=-2$, ce qui est impossible; 

soit une courbe dont les valeurs propres de Frobenius seraient les oppos\'ees
des pr\'ec\'edentes; sur $\mbox{\bf F}_{27}$, ce seraient
$-(1\pm\sqrt{-2})^3=5\pm \sqrt{-2}$, et le nombre de points serait
$27+1-30<0$, ce qui est encore impossible!

\vskip .5 truein
\hrule
\vskip .5 truein
ch\`ere Kristin,

Voici quelques compl\'ements sur le texte ``Modules hermitiens\dots''
que je vous avais envoy\'e en janvier.

Un certain nombre d'\'enonc\'es avaient \'et\'e laiss\'es sans
d\'emonstrations.  Je vais les reprendre:

1. Le plus important est:

(3.1) {\it Le foncteur $T:\Ab (a,q)\rightarrow \Mod(R)$ est une
\'equivalence de cat\'egories.}  

J'avais donn\'e une partie de
l'argument, \`a savoir la construction d'un foncteur appel\'e ``$S$''
qui transforme $L\in \Mod(R)$ en la vari\'et\'e ab\'elienne $A_L =
L\otimes_R E$.  Ce foncteur a des propri\'et\'es agr\'eables, et
faciles \`a d\'emontrer, par exemple
$$ \Hom(S(L),S(L'))   = \Hom_R(L,L').$$
Cela montre d\'ej\`a que $T\circ S$ est l'identit\'e, et donc que $T$
est injectif.  Pour montrer que $S\circ T=1$ (avec un certain abus de
notation!), on est ramen\'e \`a prouver l'\'enonc\'e:

(3.2) {\it Toute vari\'et\'e ab\'elienne $A\in \Ab(a,q)$ est de la forme
$A_L$ pour un $L$ convenable.}

C'est l\`a le point essentiel.  Il va r\'esulter de:

(3.6) {\it Soit $A\rightarrow B$ une isog\'enie dans $\Ab(a,q)$.  Si $B$
est de la forme $A_L$, alors il en est de m\^{e}me de $A$} (pour un
module $L'$ convenable).

(D'apr\`es (2.1), $A$ est isog\`ene \`a $E\times\cdots\times E$; on
applique alors (3.6) \`a $B=E\times\cdots\times E$, qui est de la
forme $A_L$, avec $L=R\times\cdots R$.)

Pour prouver (3.6), il nous faut contr\^{o}ler les isog\'enies de
vari\'et\'es ab\'eliennes, et c'est ici que l'hypoth\`ese
``ordinaire'' va \^{e}tre essentielle.  Il me faut rappeler des choses
connues.  Pla\c{c}ons-nous d'abord sur un corps alg\'ebriquement clos $k$ (on
prendra ensuite pour $k$ une cl\^{o}ture alg\'ebrique du corps fini
k).  Soit $A$ une vari\'et\'e ab\'elienne sur $k$, de dimension
$n$.  Si $l$ est un nombre premier $\not =$ caract. $k$, on sait ce
qu'est le $l$-i\`eme {\it module de Tate} $T_l(A)$ de $A$:  c'est la limite
projective des points de $l^m$-division de $A$.  C'est un ${\bf Z}_l$-module
libre de rang $2n$.  Ce module ``contr\^{o}le'' les isog\'enies
$A'\rightarrow A$ de degr\'e une puissance de $l$, au sens suivant:
une telle isog\'enie correspond bijectivement \`a un sous ${\bf Z}_l$-module
d'indice fini de $T_l(A)$ (\`a savoir l'image de $T_l(A')$ dans
$T_l(A)$).  Dans notre cas (k fini et $k$ cl\^{o}ture alg\'ebrique
de k), on voit en outre que les k-isog\'enies correspondent aux
sous-modules de $T_l(A)$ qui sont stables par l'action de Galois,
i.e. par le Frobenius $\pi$: puisqu'on a suppos\'e $A\in \Ab(a,q)$,
cela veut dire que le sous-module en question est stable par $R$.
Finalement, on voit que les k-isog\'enies $A'\rightarrow A$ qui sont
de degr\'e une puissance de $l$ sont classifi\'ees par les sous
$R_l$-modules de $T_l(A)$, o\`u $R_l = R\otimes {\bf Z}_l$.  Dans le cas
particulier $A=E$, on constate que $T_l(E)$ est un {\it $R_l$-module libre
de rang 1} (regarder les rangs!).  Si $A$ est de la forme $A_L$, on
constate aussi que $T_l(A)=T_l(E)\otimes_R L = T_l(E)\otimes L_l$,
o\`u le produit tensoriel est pris sur $R_l$ et $L_l$ d\'esigne
${\bf Z}_l\otimes L = R_l\otimes_R L$.  Il est alors imm\'ediat que les
$R_l$-sous-modules de $T_l(A)$ d'indice fini sont de la forme
$T_l(E)\otimes L'$ o\`u $L'$ est un sous-$R$-module de $L$ d'indice
une puissance de $l$, et l'on a alors $A'=E\otimes L'$.  Autrement
dit, {\it l'\'enonc\'e (3.6) est vrai si le degr\'e de l'isog\'enie est de
la forme $l^m$ avec $l\not= p$}.

On est donc ramen\'e \`a regarder le cas o\`u le degr\'e est une
puissance de la caract\'eristique $p$,  Evidemment, tout revient ici
encore \`a contr\^{o}ler les isog\'enies.  Pour une vari\'et\'e
ab\'elienne quelconque, cela peut se faire au moyen d'un {\it module de
Dieudonn\'e} convenable.  Heureusement, les vari\'et\'es ab\'eliennes
qui nous int\'eressent ici sont {\it ordinaires}, et cela simplifie
beaucoup la situation.  En effet, pour une telle vari\'et\'e $A$, on
peut d\'efinir un ``$p$-i\`eme module de Tate'' $T_p(A)$ qui a exactement
les m\^{e}mes propri\'et\'es que les $T_l$, \`a savoir: c'est un
${\bf Z}_p$-module libre de rang $2n=2.\dim(A)$, et il contr\^{o}le les
$p$-isog\'enies comme ci-dessus.  Ce module est somme directe de deux
modules de rang $n$:
$$ T_p(A) = T_p(A)_e \oplus T_p(A)_i$$
($e$= \'etale; $i$=infinit\'esimal).

La partie \'etale $T_p(A)_e$ est d\'efinie comme la limite projective
des points de $p^m$-division de $A$; la partie infinit\'esimale
$T_p(A)_i$ peut se d\'efinir comme le Hom (dans la cat\'egorie des
groupes formels) de ${\bf G}_m$ dans $A$, ou bien comme le ${\bf Z}_p$-dual de
$T_p(A^*)_e$, o\`u $A^*$ est la duale de $A$.  On d\'emontre que l'on
a les m\^{e}mes propri\'et\'es que ci-dessus pour les $T_l$.  Noter
que l'anneau $R_p = {\bf Z}_p\otimes R$ est \'egal \`a 
${\bf Z}_p\times {\bf Z}_p$ (les
deux facteurs \'etant caract\'eris\'es par le fait que $\pi$ donne une
unit\'e dans le premier, et un \'el\'ement de l'id\'eal maximal dans le
second).  Cette d\'ecomposition de $R_p$ est compatible avec la
d\'ecomposition en deux morceaux de $T_p(A)$, ainsi qu'avec le fait que
toute $p$-isog\'enie se d\'ecompose en une isog\'enie \'etale et une
isog\'enie radicielle.  Bref, tout marche tr\`es bien, et l'on arrive
ainsi \`a d\'emontrer (3.6) pour les $p$-isog\'enies.

(Une autre fa\c con de justifier (3.2) et (3.6) consiste \`a utiliser
un r\'esultat de Deligne (Invent.math. 8 (1969), 238-243) qui donne
une \'equivalence de la cat\'egorie des vari\'et\'es ab\'eliennes
ordinaires sur $\overline{{\mbox{\bf F}}}_p$ avec la cat\'egorie des
${\bf Z}$-modules libres de type fini munis d'un endomorphisme $F$ ayant un
certain nombre de propri\'et\'es raisonnables.

Cet article de Deligne contient quelques r\'ef\'erences, mais pas
beaucoup.  La th\'eorie du $T_p$ des vari\'et\'es ordinaires est
connue depuis longtemps (voir par exemple mon article de l'Amer.J. 80
(1958), 715--739), mais n'int\'eresse pas les sp\'ecialistes, car elle
est trop simple!  Le cas amusant est le cas supersingulier, \'etudi\'e
en d\'etail par Oort: on y trouve des familles ``continues'' de
$p$-isog\'enies.)

2. Il faut parler un peu de la dualit\'e et des isog\'enies ($\S\S$
4,5) L'\'enonc\'e (4.1) ne pr\'esente pas de difficult\'es.  Le point
essentiel est que, pour la courbe elliptique $E$, le transpos\'e $f^*$
d'un endomorphisme $f$ de $E$ est \'egal au {\it conjugu\'e} $\overline{f}$
de $f$, lorsqu'on identifie $f$ \`a un \'el\'ement de $R$, qui est une
extension quadratique de ${\bf Z}$.

Passons aux {\it polarisations}.  Soit $\varphi:A\rightarrow A^*$ un
morphisme, et supposons que $A\in\Ab(a,q)$ soit associ\'e au module
$L$, auquel cas $A^*$ est associ\'e \`a l'anti-dual $L^*$ de $L$.
Alors $\varphi$ correspond \`a une application $R$-lin\'eaire
$h:L\rightarrow L^*$, ou, ce qui revient au m\^{e}me, \`a une
application sesquilin\'eaire $H:L\times L\rightarrow R$.

Si $\varphi$ est une polarisation, on a $\varphi^*=\varphi$, ce qui se
traduit par $h^*=h$, ou encore par le fait que $H$ est une {\it forme
hermitienne}.

Inversement, si $H$ est une telle forme, il lui correspond
$\varphi:A\rightarrow A^*$ avec $\varphi^*=\varphi$.  Il en r\'esulte
(cf. Mumford, \textit{Abelian Varieties}, p. 188, th. 2, et p. 189,
Remarque) que $\varphi$ est de la forme ``$\varphi_D$'' pour un
diviseur $D$ sur $A$.  Dire que $\varphi$ est une polarisation
\'equivaut \`a dire que $D$ est {\it ample}.  Pour prouver (5.1) et (5.2) je
dois montrer que {\it cela se produit si et seulement si $H$ est d\'efinie
$>0$}.

Je peux remplacer $A$ par une vari\'et\'e isog\`ene.  En effet, si
$A'\rightarrow A$ est une isog\`enie, le morphisme
$\varphi:A\rightarrow A^*$ d\'efinit $\varphi':A'\rightarrow {A'}^*$ en
composant:  $A'\rightarrow A \rightarrow A^* \rightarrow {A'}^*$, et il
est facile de voir que $\varphi$ est une polarisation si et seulement
si $\varphi'$ en est une.  M\^{e}me invariance pour le fait que $H$
soit d\'efinie $>0$.  Ceci permet de choisir pour $A'$ le produit
$E\times\cdots E=E^n$, auquel cas son dual est aussi $E^n$, et
$\varphi$ est donn\'ee par une matrice hermitienne $(a_{ij})$ \`a
coefficients dans $R$.

On applique alors un r\'esultat qui se trouve dans Mumford (p.210,
lignes 4 \`a 6).  On peut aussi raisonner directement, en utilisant
encore une autre isog\'enie pour se ramener au cas o\`u la matrice
$(a_{ij})$ est une matrice {\it diagonale} avec des entiers $(d_i)$ sur la
diagonale.  Le fait que cette matrice donne une polarisation de $E^n$
si (et seulement si) les $d_i$ sont $>0$ est imm\'ediat.

Bien \`a vous.

J-P. Serre


\end{document}